# PATHWISE INEQUALITIES FOR LOCAL TIME: APPLICATIONS TO SKOROKHOD EMBEDDINGS AND OPTIMAL STOPPING

BY A. M. G. COX,[1] DAVID HOBSON[2] AND JAN OBŁÓJ[3]

*University of Bath, University of Warwick and Imperial College London*

We develop a class of pathwise inequalities of the form $H(B_t) \geq M_t + F(L_t)$, where $B_t$ is Brownian motion, $L_t$ its local time at zero and $M_t$ a local martingale. The concrete nature of the representation makes the inequality useful for a variety of applications. In this work, we use the inequalities to derive constructions and optimality results of Vallois' Skorokhod embeddings. We discuss their financial interpretation in the context of robust pricing and hedging of options written on the local time. In the final part of the paper we use the inequalities to solve a class of optimal stopping problems of the form $\sup_\tau \mathbb{E}[F(L_\tau) - \int_0^\tau \beta(B_s)\,ds]$. The solution is given via a minimal solution to a system of differential equations and thus resembles the maximality principle described by Peskir. Throughout, the emphasis is placed on the novelty and simplicity of the techniques.

**1. Introduction.** The aim of this paper is to develop and explore a new approach to solving Skorokhod embeddings and related problems in stochastic control based on pathwise inequalities of the form

$$(1) \qquad H(B_t) \geq M_t + F(L_t) \qquad \forall t \geq 0,$$

where $B_t$ is Brownian motion, $L_t$ is its local time in zero and $M_t$ is a local martingale. Then, provided the stopping time $\tau$ is finite almost surely, and provided the stopped martingale $M_{t \wedge \tau}$ is uniformly integrable, we have

$$(2) \qquad \mathbb{E}[H(B_\tau)] \geq \mathbb{E}[F(L_\tau)].$$

There is equality in (2) if there is equality at $\tau$ in (1).

Received February 2007; revised October 2007.
[1]Supported by the Nuffield Foundation.
[2]Supported by an EPSRC Advanced Fellowship.
[3]Supported by a Marie Curie Intra-European Fellowship within the 6th European Community Framework Programme.

Our aim is to find pairs $(H, F)$ such that (1) holds and to use this pathwise inequality to deduce inequalities of the form (2). We can then investigate the optimality properties of (2). For the examples we have in mind $F$ and $H$ are typically convex. Further we often consider stopping rules of the form

$$\tau_\phi = \inf\{u : B_u \notin (\phi_-(L_u), \phi_+(L_u))\}$$

and then there is a 1–1 correspondence between $\phi_{+/-}$ and the law of the stopped process.

We shall consider three different approaches to the inequalities in (1) and (2).

First, given $H$ and $\phi_{+/-}$ we find $F$ (and $M$) such that (1) and then (2) holds. This will build our intuition for constructing inequalities of this type.

Second, and more importantly, given $F$ we find $H$ such that (1) holds, and then for all $\tau$ satisfying suitable integrability conditions we also have $\mathbb{E}[F(L_\tau)] \leq \mathbb{E}[H(B_\tau)]$. If we let $\mathcal{T}(\mu)$ denote the set of stopping times such that $(B_{t \wedge \tau})$ is a uniformly integrable martingale, and such that $B_\tau \sim \mu$, then for all $\tau \in \mathcal{T}(\mu)$

$$\mathbb{E}[F(L_\tau)] \leq \int_\mathbb{R} H(x) \mu(dx).$$

In particular, for all minimal solutions of the Skorokhod embedding problem for $\mu$ in $B$ we have a bound on $\mathbb{E}[F(L_\tau)]$. We carry out this program in Section 2. We recover results of Vallois [22, 23] for Skorokhod embeddings based on local times. See Cox and Hobson [8] for a recent study concerned with similar embeddings and Obłój [15] for an extensive survey and history of the Skorokhod embedding problem.

Third, given $F$ and $H$ satisfying (1), then for suitable $\tau$ we have $\mathbb{E}[F(L_\tau) - H(B_\tau)] \leq 0$. This means we can consider problems of the form

$$\sup_\tau \mathbb{E}[F(L_\tau) - H(B_\tau)]$$

both for general $\tau$ and for $\tau \in \mathcal{T}(\mu)$ for given $\mu$; further, under suitable integrability conditions, the problem can be recast (via Itô's lemma) as the more natural stopping problem

$$\sup_\tau \mathbb{E}\left[F(L_\tau) - \int_0^\tau \tfrac{1}{2} H''(B_s)\, ds\right].$$

This is the subject of Section 3. Similar problems, but with the local time replaced by the maximum process, have been studied by Jacka [13], Dubins, Shepp and Shiryaev [9], Peskir [19], Obłój [17] and Hobson [12]. The formulation of our solution will be similar to the *maximality principle* of Peskir [19].

One of our motivations for studying inequalities of the form (2) and the relationship to pathwise inequalities such as (1), is the interpretation of such inequalities in mathematical finance as superreplication strategies for exotic derivatives, with associated price bounds. The idea is that if we can identify a martingale stock price process $S_t$ with a time-changed Brownian motion such that $S_T \sim B_\tau$, and if we know the prices of vanilla call options on $S_T$, then this is equivalent to knowing the law of $B_\tau$. If we can also identify the martingale $M$ in (1) with the gains from trade from a simple strategy in $S$, then we have a superreplicating strategy for an exotic option with payoff which is a function of the local time of $S$. Furthermore this strategy and associated price do not depend on any model assumptions.

For the case where the exotic option has a payoff which depends on the maximum (e.g., lookback and barrier options) this idea was exploited by Hobson [11] and Brown, Hobson and Rogers [5]. Financial options with payoff contingent on the local time are rare, but they can appear naturally when considering the "naive" hedging of plain vanilla options and have recently been the subject of a study by Carr [6]. A further discussion of the application of our ideas to mathematical finance is given in Section 2.3.

NOTATION. We work on a filtered probability space satisfying the usual hypotheses. $(B_t)$ denotes a real-valued Brownian motion and $(L_t)$ is its local time at zero (see Revuz and Yor [20], Chapter VI, for definition and further properties). We stress, however, that one can equally assume that $(B_t)$ is a diffusion on natural scale (thus a Markov local martingale) with $B_0 = 0$ and $\langle B \rangle_\infty = \infty$ a.s. No changes in the paper are needed apart from replacing $dt$ with $d\langle B \rangle_t$ where appropriate. Furthermore, we could also work with a recurrent time-homogenous diffusion using the scale function to change the coordinates (see the remarks in Section 3).

$F, H$ will typically denote convex functions and $\mu$ a probability measure with $\overline{\mu}(x) = \mu([x, \infty))$ denoting the right-tail. We write $X \sim \mu$ or $\mathcal{L}(X) = \mu$ to say that the law of $X$ is $\mu$.

**2. Convex functions of the terminal local time.** We begin by studying the first and second problems suggested in the Introduction. First, given $H$ and $\mu$ we will find $F$ and $M$ such that (1) holds. Then we will reverse the process, so that for a given $F$ and $\mu$ we will find $H$. For an appropriate $H$ it will follow that (2) holds for all minimal $\tau$ with $B_\tau \sim \mu$. Furthermore, the function $H$ will be optimal in the sense that there exists a stopping time $\tau_\phi$ (which we give explicitly) for which there is equality in (2). We will first consider the well-behaved case to build up intuition and then, in Section 2.2, develop the general approach.

Throughout this section we work with convex functions $F$ and $H$. The essential property of convex functions that we use is one of the most fundamental, namely that the graph of $H$ lies above any tangent to it, that is, $H(b) \geq H'(a)(b-a) + H(a)$.

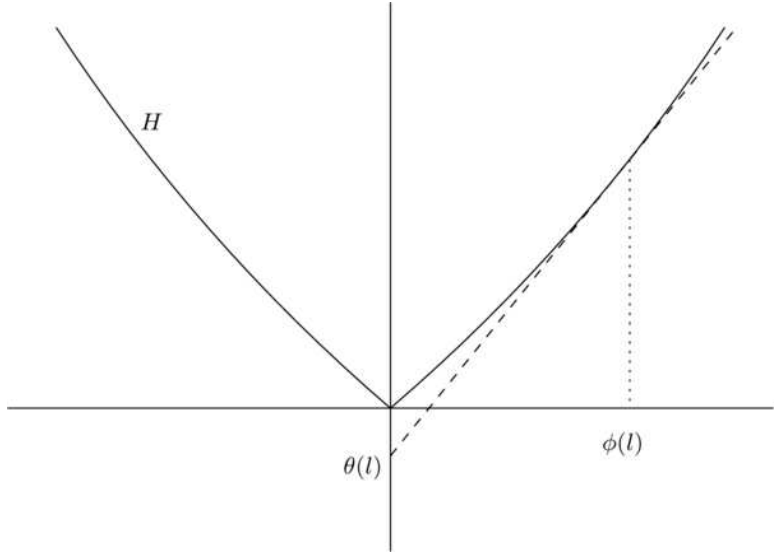

FIG. 1. *Since $H$ is convex we have $H(b) \geq H(\phi(l)) + (b - \phi(l))H'(\phi(l))$. Further the intercept of this tangent with the y-axis is $\theta(l) = H(\phi(l)) - \phi(l)H'(\phi(l))$.*

2.1. *Symmetric terminal laws with positive densities.* Let $H$ be a symmetric, strictly convex function which is differentiable on $\mathbb{R} \setminus \{0\}$. Then, for any $a, b > 0$ we have $H(-b) = H(b) \geq H'(a)(b - a) + H(a)$, with equality if and only if $b = a$; see Figure 1.

Let $\phi$ be any continuous, strictly increasing function with $\phi(0) = 0$, and let $\psi$ denote its inverse. Define $\gamma(l) = H'(\phi(l))$, $\Gamma(l) = \int_0^l \gamma(m)\,dm$ and $\theta(l) = H(\phi(l)) - \phi(l)H'(\phi(l))$. Then with $b = |B_t|$ and $a = \phi(L_t)$ we have for $t \geq 0$,

$$H(B_t) \geq H'(\phi(L_t))|B_t| - H'(\phi(L_t))\phi(L_t) + H(\phi(L_t)) \tag{3}$$
$$= M_t + F(L_t),$$

where $M_t = M_t^{H,\phi} = |B_t|\gamma(L_t) - \Gamma(L_t)$ and $F(l) = F_{H,\phi}(l) = \Gamma(l) + \theta(l)$.

By construction $M_t^{H,\phi}$ is a local martingale (cf. Obłój [16]), so if $\tau$ is a stopping time such that $\mathbb{E}[M_\tau^{H,\phi}] = 0$, then

$$\mathbb{E}[F_{H,\phi}(L_\tau)] \leq \mathbb{E}[H(B_\tau)]. \tag{4}$$

Define $\tau_\phi = \inf\{u > 0 : |B_u| = \phi(L_u)\}$ and suppose $\phi$ is such that $0 < \tau_\phi < \infty$ a.s. Let $\mu = \mathcal{L}(B_{\tau_\phi})$. Then for any solution of the Skorokhod embedding problem for $\mu$ in $B$ with the property that $\mathbb{E}[M_\tau^{H,\phi}] = 0$ we have

$$\mathbb{E}[F_{H,\phi}(L_\tau)] \leq \int_\mathbb{R} H(x)\mu(dx). \tag{5}$$

Thus we obtain an upper bound for the value $\mathbb{E}[F_{H,\phi}(L_\tau)]$. There is equality in (5) if $\tau = \tau_\phi$ and $\mathbb{E}[M^{H,\phi}_{\tau_\phi}] = 0$. Further, among $\tau \in \mathcal{T}(\mu)$ with $\mathbb{E}[M^{H,\phi}_\tau] = 0$, this is the only stopping time with this property. To see this recall that a stopping time $\tau$ is minimal if and only if $(B_{t \wedge \tau})$ is uniformly integrable (Monroe [14], Theorem 3), and since we have strict inequality in (3) unless $|B_\tau| = \phi(L_\tau)$, it must be the case that for $\tau \in \mathcal{T}(\mu)$ to yield equality in (5) we must have that $\tau$ is the first positive time that $|B_\tau| = \phi(L_\tau)$.

Stopping times of the form $\tau_\phi$ were used by Vallois [22] to solve the Skorokhod embedding problem. For a given symmetric centered probability measure $\mu$ on $\mathbb{R}$, Vallois ([22], Theorem 3.1) defined a function $\phi = \phi(\mu)$ such that $B_{\tau_\phi} \sim \mu$ and $(B_{t \wedge \tau_\phi} : t \geq 0)$ is a uniformly integrable martingale.[4] Vallois [23], Theorem 1, then proved that his stopping times maximize the expectation of convex functions of $L_T$. Using our methodology we recover his results: both the embedding property and the optimality with respect to the convex ordering. Formally, the results of this section are not new but the emphasis is on the novelty of the method. Indeed, in both of his papers Vallois relied on martingale methods to compute laws of stopped processes and did not have any pathwise inequalities in the spirit of (3). In this sense our method offers a new interpretation and novel applications, in particular in the context of financial mathematics (see Section 2.3 below).

Our aim now is to reverse the procedure described above. For a given convex function $F$ and measure $\mu$ we aim to find $H$ such that $F_{H,\phi} = F$ where $\phi$ is related to $\mu$ via Vallois' solution to the Skorokhod embedding problem.

To illustrate our method we begin with the simplest case and for the remainder of this section we adopt the following simplifying assumptions:

(A1.1) Suppose $\mu$ is symmetric, with finite first moment and with a positive density (with respect to Lebesgue measure) on $\mathbb{R}$.
(A1.2) Suppose $F : \mathbb{R}^+ \mapsto \mathbb{R}^+$ is convex and increasing and has continuous first derivative $F'$ which is bounded by $K$.
(A1.3) Suppose $\phi$ is any continuous, strictly increasing function such that $\phi(0) = 0$ and such that $\int_{0+} dl/\phi(l) < \infty$ and $\int^\infty dl/\phi(l) = \infty$.

Denote the inverse to $\phi$ by $\psi$. Define the measure $\nu \equiv \nu_\phi$ on $\mathbb{R}^+$ via

$$(6) \qquad \overline{\nu}(l) = \nu([l, \infty)) = \exp\left(-\int_0^l \frac{dm}{\phi(m)}\right).$$

By the assumptions on $\phi$, $\nu$ is a probability measure with density $\overline{\nu}(l)/\phi(l)$.

---
[4]Vallois [22] actually solved the problem for any centered probability measure on $\mathbb{R}$ considering more general asymmetric stopping times as in (16) below. See [22] or Obłój ([15], Section 3.12) for more details on the general construction.

Define $H \equiv H_{F,\phi}$ via

$$H'(b) = \frac{1}{\overline{\nu}(\psi(b))} \int_{\psi(b)}^{\infty} F'(m)\nu(dm); \qquad H(0) = F(0). \tag{7}$$

By the assumptions on $F$ we have that $H'$ is well defined and bounded by $K$. It is easy to show that $H'$ is increasing so that $H$ is convex. Indeed,

$$H''(b) = \frac{\psi'(b)}{b\overline{\nu}(\psi(b))} \int_{\psi(b)}^{\infty} [F'(m) - F'(\psi(b))]\nu(dm) \tag{8}$$

$$= \frac{\psi'(b)}{b}(H'(b) - F'(\psi(b))) \tag{9}$$

which is nonnegative since $\psi$ is increasing and $F$ is convex.

LEMMA 2.1. *Suppose $H'(x) \leq K$ and $\tau$ is such that $(B_{t \wedge \tau})$ is a uniformly integrable (UI) martingale. Then $\mathbb{E}[M_\tau^{H,\phi}] = 0$.*

PROOF. Let $\sigma_n = \inf\{t : |B_t| \geq n\}$, $\rho_m = \inf\{t : L_t \geq m\}$ and $\tau_{m,n} = \tau \wedge \sigma_n \wedge \rho_m$. As the local martingale $(M_{t \wedge \tau_{m,n}}^{H,\phi} : t \geq 0)$ is bounded it is UI and $\mathbb{E} M_{\tau_{m,n}}^{H,\phi} = 0$ so that

$$\mathbb{E}\Gamma(L_{\tau_{m,n}}) = \mathbb{E}\gamma(L_{\tau_{m,n}})|B_{\tau_{m,n}}| = \mathbb{E}\gamma(L_{\tau_{\infty,n}})|B_{\tau_{\infty,n}}|\mathbf{1}_{\tau \wedge \sigma_n \leq \rho_m}.$$

By the monotone convergence theorem, as $m \to \infty$ both sides converge, and in the limit we obtain $\mathbb{E}\Gamma(L_{\tau_{\infty,n}}) = \mathbb{E}\gamma(L_{\tau_{\infty,n}})|B_{\tau_{\infty,n}}|$. Now, as $n \to \infty$, the left-hand side converges, again by the monotone convergence theorem, to $\mathbb{E}[\Gamma(L_\tau)]$, since $\tau = \tau_{\infty,\infty}$. The right-hand side converges to $\mathbb{E}\gamma(L_\tau)|B_\tau|$ since $\gamma$ is bounded and $|B_{t \wedge \tau}|$ is UI, so that finally $\mathbb{E}[M_\tau^{H,\phi}] = 0$. □

PROPOSITION 2.2. (i) *Define $H \equiv H_{F,\phi}$ via (7). Then, for all $\tau$ such that $(B_{t \wedge \tau})$ is a uniformly integrable martingale,*

$$\mathbb{E}[F(L_\tau)] \leq \mathbb{E}[H_{F,\phi}(B_\tau)]. \tag{10}$$

(ii) *Let $\phi_\mu$ be the inverse to $\psi_\mu$ where $\psi_\mu$ is given by*

$$\psi_\mu(x) = \int_0^x \frac{s}{\overline{\mu}(s)} \mu(ds). \tag{11}$$

*Let $\tau_\mu \equiv \tau_{\phi_\mu} = \inf\{u > 0 : |B_u| = \phi_\mu(L_u)\}$. Then $B_{\tau_\mu} \sim \mu$, and $\mathbb{E}[F(L_{\tau_\mu})] = \int H_{F,\phi_\mu}(x)\mu(dx)$.*

(iii) $\forall \tau \in \mathcal{T}(\mu)$, $\mathbb{E}[F(L_\tau)] \leq \mathbb{E}[F(L_{\tau_\mu})]$.

PROOF. (i) The first part follows from (4) provided we can show that $F_{H,\phi} \equiv F$ and $\mathbb{E}[M_\tau^{H,\phi}] = 0$. This latter statement is guaranteed by Lemma 2.1. For the former, recall that $F_{H,\phi}(l) = \int_0^l H'(\phi(m))\,dm - \phi(l)H'(\phi(l)) + H(\phi(l))$. Setting $l = \psi(b)$ and differentiating, we obtain from (9)

$$\psi'(b)F'_{H,\phi}(\psi(b)) = \psi'(b)H'(b) - bH''(b) = \psi'(b)F'(\psi(b)).$$

Since $F_{H,\phi}(0) = H(0) = F(0)$ and the image of $\psi$ is the whole of $\mathbb{R}$ we conclude that $F_{H,\phi} \equiv F$.

(ii) Note first that

$$\int_0^u \frac{dl}{\phi_\mu(l)} = \int_0^{\phi_\mu(u)} \frac{\mu(ds)}{\overline{\mu}(s)} = -\log(\overline{\mu}(\phi_\mu(u))),$$

which is finite for $u \in (0, \infty)$ and infinite for $u = \infty$. Hence $\phi_\mu$ satisfies Assumption (A1.3).

Now let $\phi$ be any function which satisfies Assumption (A1.3), and let $\tau_\phi = \inf\{t > 0 : |B_t| \geq \phi(L_t)\}$. By an excursion theory argument (cf. Obłój and Yor [18])

(12) $$\mathbb{P}(L_{\tau_\phi} > l) = \exp\left(-\int_0^l \frac{ds}{\phi(s)}\right), \qquad l > 0,$$

and, using Assumption (A1.3), we have $0 < L_{\tau_\phi} < \infty$ a.s. and therefore also $0 < \tau_\phi < \infty$ a.s. We have $|B_{\tau_\phi}| = \phi(L_{\tau_\phi})$ and, as remarked earlier, equality is achieved in (3), so that $\mathbb{E}F(L_{\tau_\phi}) = \mathbb{E}H(B_{\tau_\phi})$.

It remains to show that for the choice $\phi = \phi_\mu$, the law of $B_{\tau_\mu}$ is $\mu$. Write $\rho$ for the law of $|B_{\tau_\mu}|$. To see directly that $\rho = 2\mu_{|\mathbb{R}_+}$ write $\overline{\rho}(x) = \mathbb{P}(L_{\tau_\mu} \geq \psi_\mu(x))$ which can be computed via (12) and (11). We want to give, however, a natural approach to recover (11) where we only suppose that $\mathbb{E}|B_{\tau_\mu}| < \infty$. We know from the comments about equality in (5) that for a wide class of functions $H$,

(13) $$\int_0^\infty F_{H,\phi}(\psi(x))\rho(dx) = \mathbb{E}F_{H,\phi}(L_{\tau_\mu}) = \mathbb{E}H(|B_{\tau_\mu}|) = \int_0^\infty H(x)\rho(dx).$$

This holds in particular for $H(x) = (|x| - k)^+$ and then the right-hand side is finite. We have

$$F_{H,\phi}(\psi(x)) = \int_0^{\psi(x)} H'(\phi(u))\,du + H(x) - xH'(x)$$
$$= \int_0^x H'(u)\,d\psi(u) + H(x) - xH'(x),$$

which substituted into (13) yields

$$\int_0^\infty \rho(dx) \int_0^x H'(u)\,d\psi(u) = \int_0^\infty xH'(x)\rho(dx).$$

Changing the order of integration we conclude

$$\int_0^\infty \overline{\rho}(x) H'(x)\, d\psi(x) = \int_0^\infty x H'(x) \rho(dx).$$

Given that the family of functions $H'(x)x$ contains the functions $f_k(x) = x\mathbf{1}_{x \geq k}$, for all $k \geq 0$, and that this family is rich enough to determine probability measures on $\mathbb{R}_+$, it follows that

$$\frac{d\psi(x)}{x} = \frac{\rho(dx)}{\overline{\rho}(x)}.$$

In particular, if $\psi \equiv \psi_\mu$ so that $\psi$ solves (11), then $d(\log(\overline{\rho}(x))) = d(\log(\overline{\mu}(x)))$ and thus $\overline{\rho}(x) = 2\overline{\mu}(x)$ where the constant 2 arises from the fact that $\overline{\rho}(0) = 1 = 2\overline{\mu}(0)$. Since $\mathcal{L}(B_{\tau_\phi})$ is symmetric we conclude $B_{\tau_\mu} \sim \mu$.

(iii) This follows immediately from (i) and (ii). □

REMARK. From the definition of $\psi_\mu$ we have for its inverse $\phi_\mu$

$$d[\ln \overline{\mu}(\phi_\mu(l))] = -\frac{dl}{\phi_\mu(l)} = d[\ln \overline{\nu}_\mu(l)],$$

where $\nu_\mu$ is defined from (6) using $\phi_\mu$. It follows that $\overline{\nu}_\mu(l) = 2\overline{\mu}(\phi_\mu(l))$. Hence (7) can be rewritten as

$$H'(b) = \frac{1}{\overline{\mu}(b)} \int_b^\infty F'(\psi(x))\mu(dx); \qquad H(0) = F(0).$$

COROLLARY 2.3. *Suppose the Assumption* (A1.2) *is relaxed, and we assume only that $F$ is convex and increasing. Then, for all $\tau \in \mathcal{T}(\mu)$, $\mathbb{E}[F(L_\tau)] \leq \mathbb{E}[F(L_{\tau_\phi})]$. In particular, the assumptions that $F'$ is continuous and $F' \leq K$ can be removed.*

PROOF. It is clear from the definition of $H$ via (7), and the proof of convexity in (8) that we do not need the derivative $F'$ continuous, but just that the integrals in (7) and (8) are well defined. Further, for any increasing convex function $F$ we define $F_K$ via $F_K(0) = F(0)$ and $F'_K = F' \wedge K$. We have shown so far that for any solution to the Skorokhod embedding problem $\mathbb{E}F_K(L_{\tau_\phi}) \geq \mathbb{E}F_K(L_\tau)$. Taking limits as $K \to \infty$, via the monotone convergence theorem, we obtain the general optimal property of Vallois' stopping time: $\mathbb{E}F(L_{\tau_\phi}) \geq \mathbb{E}F(L_\tau)$ for smooth symmetric terminal distributions $B_\tau \sim B_{\tau_\phi} \sim \mu$. □

EXAMPLE 2.4. Suppose $\overline{\mu}(l) = e^{-2\alpha^2 l}/2$, where $\alpha > 0$. Then $\psi_\mu(b) = \alpha^2 b^2$, $\phi_\mu(l) = \sqrt{l}/\alpha$ and $\overline{\nu}(l) = 2\overline{\mu}(\phi(l)) = e^{-2\alpha\sqrt{l}}$.

Suppose now that $H(x) = Ax^2 + B|x|$, with $A, B \geq 0$. Then

$$F_{H,\phi}(l) = \frac{4}{3\alpha} A l^{3/2} + \left(B - \frac{1}{\alpha^2}A\right)l.$$

Note that $F$ is convex and increasing if and only if $B \geq A/\alpha^2$. Conversely, if $F(l) = Cl^{3/2}$, then $H_{F,\phi}(x) = (3\alpha/4)Cx^2 + (3/4\alpha)C|x|$.

2.2. *Arbitrary centered terminal laws.* We want to extend the results of the previous section to arbitrary terminal laws. We need to be able to deal with two issues: atoms in $\mu$ and asymmetry of $\mu$. We deal with the former by parameterizing the fundamental quantities in terms of the quantiles of $\mu$ and we deal with the latter by introducing separate functions $\phi_+$ and $\phi_-$ on the positive and negative half-spaces, respectively.

*In this section when we take inverse functions we always mean the right-continuous versions. We also use the notation $\phi_\pm$ to indicate the pair $(\phi_+, \phi_-)$, this should cause no confusion.*

Let $\phi_+ : \mathbb{R}_+ \to \mathbb{R}_+$ be an increasing function and $\phi_- : \mathbb{R}_+ \to \mathbb{R}_-$ a decreasing function. To develop an analogue to (3) we will parameterize the negative half-line with $\phi_-$ and the positive half-line with $\phi_+$ so that

$$\begin{cases} H(x) \geq \gamma_+(l)x + \theta_+(l), & x > 0, \\ H(z) \geq \gamma_-(l)z + \theta_-(l), & z < 0, \end{cases}$$

where (see Figure 2)

$$\begin{cases} H'(\phi_+(l-)) \leq \gamma_+(l) \leq H'(\phi_+(l+)), & \theta_+(l) = H(\phi_+(l)) - \phi_+(l)\gamma_+(l), \\ H'(\phi_-(l+)) \leq \gamma_-(l) \leq H'(\phi_-(l-)), & \theta_-(l) = H(\phi_-(l)) - \phi_-(l)\gamma_-(l). \end{cases}$$

Substituting $B_t$ and $L_t$ we obtain

$$\begin{aligned}
(14) \quad H(B_t) &\geq \gamma_+(L_t)B_t^+ - \gamma_-(L_t)B_t^- + \theta_+(L_t)\mathbf{1}_{B_t \geq 0} + \theta_-(L_t)\mathbf{1}_{B_t < 0} \\
&= M_t^{H,\phi} + \Gamma(L_t) + \theta_+(L_t)\mathbf{1}_{B_t \geq 0} + \theta_-(L_t)\mathbf{1}_{B_t < 0},
\end{aligned}$$

where $\Gamma(l) = \int_0^l (\gamma_+(m) - \gamma_-(m))/2 \, dm$ and $M_t^{H,\phi} = \gamma_+(L_t)B_t^+ - \gamma_-(L_t)B_t^- - \Gamma(L_t)$ is a local martingale. If we choose the various quantities such that $\theta_+(l) = \theta_-(l) = \theta(l)$, then we have

$$(15) \qquad H(B_t) \geq M_t^{H,\phi} + F_{H,\phi,\gamma}(L_t),$$

where $F_{H,\phi,\gamma}(l) = \Gamma(l) + \theta(l)$.

Note that when $H'$ is not continuous different choices of functions $\gamma_\pm(l)$ (or equivalently different choices of tangents to $H$) may lead to different functions $F_{H,\phi,\gamma}(l)$ and different inequalities.

As before, our goal is to reverse the procedure: given $F$ and a centered probability measure $\mu$ on $\mathbb{R}$, we aim to choose $H$, $\phi$ and $\gamma$ such that $F_{H,\phi,\gamma} \equiv F$ and $B_{\tau_\phi} \sim \mu$ where

$$(16) \qquad \tau_\phi = \inf\{u > 0 : B_u^+ = \phi_+(L_u) \text{ or } B_u^- = -\phi_-(L_u)\}.$$

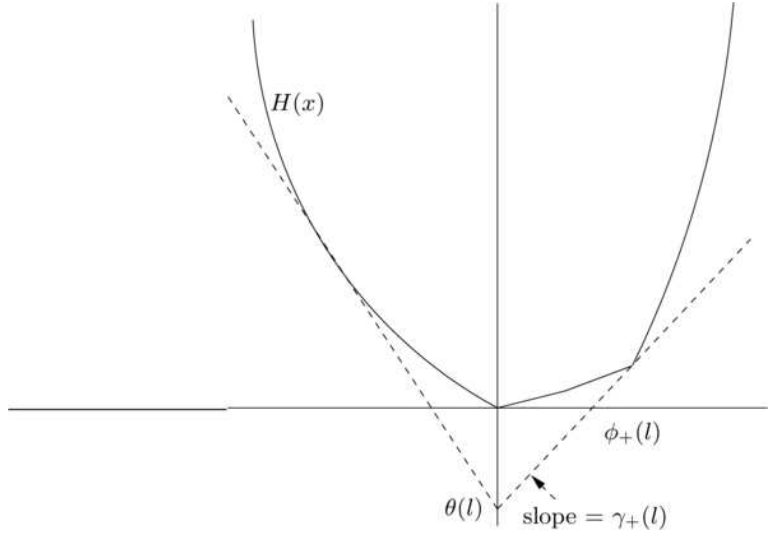

Fig. 2. *Specification of the various functions, $H$, $\phi_+$, $\gamma$ and $\theta$.*

Define

$$\Delta(l) = \int_0^l \left( \frac{1}{2\phi_+(m)} + \frac{1}{2|\phi_-(m)|} \right) dm.$$

ASSUMPTION (A2). Suppose $F$ is convex and increasing and suppose $\phi_+$ and $\phi_-$ are increasing positive and decreasing negative functions, respectively, such that $\Delta(l)$ is finite for each $l > 0$, but increases to infinity with $l$.

Fix $\phi_+$ and $\phi_-$ satisfying Assumption (A.2), and let $\psi_+, \psi_-$ denote their respective inverses. Define $\nu = \nu_\phi$ via $\overline{\nu}(l) = \exp(-\Delta(l))$. Given $F$, define the increasing function $\Sigma$ via

$$\Sigma(l) = \frac{1}{\overline{\nu}(l)} \int_l^\infty F'(m) \nu(dm).$$

Where $F'(l)$ exists, define $\delta(l) = \Sigma(l) - F'(l)$. Then $\delta$ is defined almost everywhere in $l$ and is positive.

Define

$$A_+(l) = \Sigma(0) + \int_0^l \frac{\delta(m)}{\phi_+(m)} dm, \qquad C(l) = -F(0) + \int_0^l \delta(m) \, dm,$$

$$A_-(l) = -\Sigma(0) - \int_0^l \frac{\delta(m)}{|\phi_-(m)|} dm,$$

and the function $H$ via

$$(17) \quad H(x) = \begin{cases} \sup_{l>0}\{xA_+(l) - C(l)\}, & x > 0, \\ F(0), & x = 0, \\ \sup_{l>0}\{xA_-(l) - C(l)\}, & x < 0. \end{cases}$$

REMARK. In fact, the only condition we need on $A_+(0)$ and $A_-(0)$ is that $A_+(0) - A_-(0) = 2\Sigma(0)$ and $A_+(0)$ and $A_-(0)$ are undetermined except through this difference. However, we fix both of them using an antisymmetry condition. A different convention for the choice of $A_+(0)$ would lead to a modification $H(x) \mapsto H(x) + xk$ for some constant $k$. For $\tau \in \mathcal{T}(\mu)$ we have $k\mathbb{E}[B_\tau] = 0$, and hence such a modification would have no effect on inequalities such as (2).

LEMMA 2.5. *$H$ is convex. The suprema for $x > 0$ and $x < 0$ in* (17) *are attained at $l = \psi_+(x)$ and $l = \psi_-(x)$, respectively. Further $H'(\phi_+(l-)) \le A_+(l) \le H'(\phi_+(l+))$ and $H'(\phi_-(l+)) \le A_-(l) \le H'(\phi_-(l-))$.*
*Finally, $F_{H,\phi,A} = F$.*

PROOF. We have

$$xA_+(l) - C(l) = x\Sigma(0) + F(0) + \int_0^l \delta(m)\left(\frac{x}{\phi_+(m)} - 1\right)dm$$

which is maximized by $l = \psi_+(x)$ since thereafter the integrand is negative. Convexity of $H$ follows immediately from the definition as a supremum of linear functions (cf. Hiriart-Urruty and Lemaréchal [10], Section B.2.1). Note also that $H$ is continuous at 0.

For the final statement observe that by definition $F_{H,\phi,A}(l) \equiv \Gamma_A(l) + \theta(l)$. We have $\theta_+(l) = \theta_-(l) = -C(l)$ and

$$A_+(l) - A_-(l) = 2\Sigma(0) + 2\int_0^l \delta(m)\frac{\nu(dm)}{\overline{\nu}(m)}$$

$$= 2\Sigma(0) + 2\int_0^l \frac{\nu(dm)}{\overline{\nu}(m)^2}\int_m^\infty [F'(m) - F'(l)]\nu(dm)$$

$$= 2\Sigma(0) + 2\int_0^l \Sigma'(m)\,dm = 2\Sigma(l).$$

As a consequence, $F_{H,\phi,A}(l) = \int_0^l \Sigma(m)\,dm - C(l) = F(0) + \int_0^l F'(m)\,dm = F(l)$. □

We can now deduce our theorem which makes precise the ideas outlined in the Introduction.

THEOREM 2.6. *Suppose $F$ and $\phi$ satisfy Assumption* (A2). *Define* $H \equiv H_{F,\phi,A}$ *via* (17). *Then, for all $\tau$ such that $(B_{t\wedge\tau})$ is a uniformly integrable martingale*

(18) $$\mathbb{E}[F(L_\tau)] \leq \mathbb{E}[H(B_\tau)].$$

PROOF. Suppose $F' \leq K$ (the result for the general case can be deduced as in Corollary 2.3). Then, by a slight generalization of Lemma 2.1, $\mathbb{E}[M_\tau^{H,\phi}] = 0$. The result now follows from (15). □

Our goal is to prove that there can be equality in (18). Moreover, if given a centered distribution $\mu$ we can find a stopping rule such that $B_\tau \sim \mu$ and there is equality in (18), then, as in Proposition 2.2, we have a tight bound on $\mathbb{E}[F(L_\tau)]$ over solutions of the Skorokhod embedding problem for $\mu$. The existence and form of an embedding of $\mu$ based on the local time, and its optimality in the sense of maximizing convex functions, are due to Vallois [22], Théorème 3.1 and [23], Théorème 1.

Let $\mu$ be a centered probability distribution with no atom in zero and let $\mu(\mathbb{R}^-) = a_* > 0$. Let $G$ denote the cumulative distribution function of $\mu$ so that $G(x) = \mu((-\infty, x])$. For $a_* \leq a \leq 1$ define $\alpha(a)$ via

$$\int_{a_*}^a G^{-1}(c)\, dc + \int_{\alpha(a)}^{a_*} G^{-1}(c)\, dc = 0.$$

Then $0 < \alpha(a) < a_*$, $\alpha(a_*) = a_*$, $\alpha(1) = 0$ and $\alpha$ is a strictly decreasing absolutely continuous function with $\alpha'(c) = G^{-1}(c)/G^{-1}(\alpha(c))$.

Define $\xi = \xi_\mu$ via

$$\xi(a) = 2\int_{a_*}^a \frac{G^{-1}(c)}{\alpha(c) + (1-c)}\, dc, \qquad a_* \leq a \leq 1,$$

and

$$\xi(a) = 2\int_a^{a_*} \frac{G^{-1}(c)}{c + 1 - \alpha^{-1}(c)}\, dc, \qquad 0 \leq a \leq a_*.$$

Then, $\xi$ is an absolutely continuous function which is convex on $a \geq a_*$ and concave on $a \leq a_*$. Note also that $\xi(1) = \infty = -\xi(0)$.

Define $\psi_\mu(x)$ via $\psi_\mu(x) = \xi(G(x))$ so that $\psi_\mu$ is an increasing function on $\mathbb{R}$. If $\mu$ is symmetric, then $\alpha(c) = 1 - c$ and for $x > 0$ we obtain the following generalization of the formula (11):

$$\psi_\mu(x) = -\int_{[0,x]} s\, d(\ln \overline{\mu}(s)).$$

THEOREM 2.7. *Let $\mu$ be a centered probability distribution on $\mathbb{R}$ with $\mu(\{0\}) = 0$. Let $\phi_\mu : \mathbb{R} \mapsto \mathbb{R}$ be the inverse to $\psi_\mu$ defined above, and define $\phi_+ : \mathbb{R}^+ \mapsto \mathbb{R}^+$ and $\phi_- : \mathbb{R}^+ \mapsto \mathbb{R}^-$ via $\phi_\pm(l) = \phi_\mu(\pm l)$. Define $\tau_\phi$ as in (16). Then $B_{\tau_\phi} \sim \mu$ and there is equality in (18).*

PROOF. First we show that $\overline{\nu}(\xi(a)) = \alpha(a) + 1 - a$ for $a > a_*$. We have

$$\ln \overline{\nu}(\xi(a)) = -\Delta(\xi(a)) = -\int_{a_*}^{a} \xi'(c) \Delta'(\xi(c))\,dc$$

$$= \int_{a_*}^{a} \frac{-1}{(\alpha(c) + 1 - c)} \left(1 - \frac{G^{-1}(c)}{G^{-1}(\alpha(c))}\right) dc$$

$$= \ln(\alpha(a) + 1 - a),$$

where we use the fact that identities $\phi_+(\xi(c)) = G^{-1}(c)$ and $\phi_-(\xi(c)) = G^{-1}(\alpha(c))$ hold $dc$ almost everywhere on $(a_*, 1)$. This implies

$$\lim_{m \uparrow \infty} \Delta(m) = \infty$$

as required by Assumption (A.2).

For $x > 0$ we have

$$\mathbb{P}(B_{\tau_\phi} > x) = \int_0^\infty \overline{\nu}(l) \frac{dl}{2\phi_+(l)} \mathbf{1}_{\{\phi_+(l) > x\}}$$

$$= \int_{(x,\infty)} \frac{\overline{\nu}(\xi(G(y)))}{2y}\, d\xi(G(y))$$

$$= \int_{G(x)}^{1} \frac{\overline{\nu}(\xi(c))}{\alpha(c) + (1-c)}\, dc = 1 - G(x).$$

Calculations for $a \in (0, a_*)$ and for $x < 0$ are similar.

Equality in (18) follows from the definition of $\tau_\phi$ and resulting equality in (14). □

COROLLARY 2.8. *Suppose $\mu$ is centered with $\mu(\{0\}) = 0$ and $F$ is convex. Then, for all $\tau \in \mathcal{T}(\mu)$, $\mathbb{E}[F(L_\tau)] \leq \mathbb{E}[F(L_{\tau_\phi})]$.*

Finally, we relax the condition that $\mu(\{0\}) = 0$. If $\mu$ places mass at zero, then we can construct an embedding of $\mu$ as follows. Let $Z$ be a Bernoulli random variable with $\mathbb{P}(Z = 0) = \mu(\{0\})$ which is independent of $B$—if necessary we expand the probability space so that it is sufficiently rich as to support $Z$—and, given $X \sim \mu$, let $\tilde{\mu}$ be the law of $X$ conditioned to be nonzero.

On $Z = 0$ set $\tau = 0$. Otherwise, on $Z = 1$, let $\tau$ be the stopping rule defined via $\tilde{\phi}$ and (16), where $\tilde{\phi}$ is defined from $\tilde{\mu}$ using the algorithm described following Theorem 2.6.

It is clear $B_\tau \sim \mu$. Also it is clear that with $H$ defined relative to $F$ and $\tilde{\phi}$, (18) still holds. Further, by considering the cases $Z = 0$, $Z = 1$ separately we see that there can be equality in (18).

Independent randomization is necessary for a stopping rule to attain equality in (18). Otherwise, if we insist that the stopping times are adapted to the minimal filtration generated by $B$, then the best that is possible is to find a sequence of times $\tau_n$ such that $B_{\tau_n} \sim \mu$ and $\lim_{n \uparrow \infty} \mathbb{E}[F(L_{\tau_n})] = \int H(x) \mu(dx)$.

2.3. *Financial applications.* Let $S_t$ be the time-$T$ forward-price process of a financial asset. (To keep notation simple we express all prices in terms of monetary units at time $T$.) Consider the following "naive" hedging strategy for a European call option with maturity $T$ and strike $K \geq S_0$: borrow $K$ and trade such that the portfolio holdings are $\max\{S_t, K\}$. In particular, the first time, if ever, that the forward exceeds $K$, buy the forward; if subsequently the forward price falls below $K$, then sell; whereupon the process is repeated. Provided $S_t$ is continuous, all the transactions happen when $S_t = K$.

Such a strategy was called the *stop-loss start-gain strategy* by Seidenverg [21]. At maturity this strategy yields $K + (S_T - K)^+$ and paying back $K$ we have replicated the call payoff at no cost. Therefore, for no arbitrage to hold the price of an out-of-the-money call would have to be zero, while in practice such calls have positive value. The answer to the apparent paradox is that when $S_t$ is continuous but has unbounded variation, then trading continuously at level $K$ accumulates local time at that level, and the strategy is not self-financing.

This resolution of the paradox, identified by Carr and Jarrow [7], shows that local time related quantities can arise naturally in financial markets. Other products closely linked with local time include corridor variance swaps, or more generally products dependent on number of downcrossings of an interval (see also Carr [6]). When exposed to the risk related to the local time, in addition to model-based prices one would want to have *model-free* bounds on the risk-quantifying products, and our study can be interpreted in this way. Analogous studies based on the supremum process (and yielding the Azéma–Yor solution to the Skorokhod embedding problem) led to *model-free* bounds on prices of look-back and barrier options (cf. Hobson [11] and Brown, Hobson and Rogers [5]).

We work in a financial market which admits *no arbitrage*, so that there exists a risk-neutral measure (equivalent to the physical measure) under which the forward-price process $(S_t)$ is a local martingale. We further assume that $(S_t)$ has continuous paths and is a true martingale (and thus a UI martingale) on finite time horizon $[0, T]$. In today's markets plain vanilla options are traded liquidly and it is an established practice to use them to calibrate models. We assume that $T$-maturity calls with the full continuum of possible strikes are traded, and then from Breeden and Litzenberger [4] and subsequent works, we know that differentiating the maturity-$T$ call prices twice with respect to the strike, we recover the probability distribution of $S_T$ under the risk-neutral measure.

To apply our results directly we define a shifted process $P_t := S_t - S_0$ with initial value $P_0 = 0$. Under the risk-neutral measure, $S_t$ is a continuous martingale and the distribution of $P_T$ is a centered distribution on $[-S_0, \infty)$ which we denote by $\mu$. As we now show, our results give bounds on the value of a contingent claim paying $F(L_T(P))$ at time $T$, where $F$ is some convex function and $(L_t(P))$ is the local time in zero of $(P_t)$, which is also the local time of $(S_t)$ at the level $S_0$.

The process $(P_t : t \leq T)$ can be written as a time-changed Brownian motion $(B_{\tau_t} : t \leq T)$ where $\tau = \tau_T$ is a stopping time such that $B_\tau \sim \mu$ and $(B_{\tau \wedge s} : s \geq 0)$ is a UI martingale. Furthermore, $L_T(P)$ is equal to the stopped Brownian local time $L_\tau$ (cf. Obłój [16]). Theorems 2.6 and 2.7 imply that $\Theta = \int H(x) \mu(dx)$, where $H$ is given explicitly via (17) for $\phi_\pm$ as in Theorem 2.7, is the upper *model-free* bound on the expected value of $F(L_T(P))$.

Associated with the price bound is a superreplicating portfolio, consisting of a static portfolio paying $H(P_T)$ and a dynamic hedge. The European payoff $H(P_T)$ can be written as a static portfolio of puts and calls. The dynamic component is given by a self-financing portfolio $G_t$ whose increase is given by $dG_t = -\Delta_t \, dS_t$ where

$$\Delta_t = H'(\phi_+(L_t(P)))\mathbf{1}_{P_t > 0} + H'(\phi_-(L_t(P)))\mathbf{1}_{P_t > 0}. \tag{19}$$

Note that $G_t$ is simply the time-change of the martingale $-M_t^{H,\phi}$. Then (15) implies that $F(L_T(P)) \leq H(P_T) + G_T$ a.s. and we have exhibited a superreplicating portfolio. The portfolio holdings $\Delta_t$ are only rebalanced when $S_t = S_0$, so that the hedging strategy is comparatively simple compared to a full dynamic hedging strategy in a Black–Scholes style model.

This approach gives an upper bound on the potential model-based prices of options contingent upon local time. The pricing mechanism in which the price of the security paying $F(L_T(P))$ is set to be $\Theta$ may be too conservative, but it does have the benefit of being associated with a superhedging strategy which is guaranteed to be successful, pathwise. A selling price $\tilde{\Theta} < \Theta$ can only be justified if the forward price process is known to belong to some subclass of models. Even in this case the seller can still use the hedging mechanism described above and be certain that his potential loss is bounded below by $\Theta - \tilde{\Theta}$ regardless of all other factors.

**3. Optimal stopping problems.** In this section we consider related optimal stopping problems. In particular, we consider solutions to problems of the form

$$\sup_\tau \mathbb{E}\left[F(L_\tau) - \int_0^\tau \beta(B_s) \, ds\right], \tag{20}$$

subject to the expectation of the integral term being finite. Both in terms of the function we wish to maximize, and the form that our solution will take, this problem can be considered as a relative of the problem considered in Peskir [19]; in particular, our solution resembles the *maximality principle* introduced by Peskir [19]. We assume (initially) only that $F$ and $\beta$ are both nonnegative; we will make stronger assumptions later as required. As stressed in the Introduction, in what follows Brownian motion $(B_t)$ could be replaced with a diffusion in natural scale $(X_t)$. We then replace $ds$ with $d\langle X \rangle_s$ where appropriate, in particular in (20), but no other changes are needed.[5] We note, however, that the continuity and time-homogeneity of the process are important. We could not easily deal, for example, with jumps (for a discussion of optimization problems for processes including jumps see, e.g., Alili and Kyprianou [1]).

Our motivation in studying (20) is threefold. First, it is a natural counterpart to the study of similar problems where the local time is replaced with the unilateral supremum process, undertaken in particular by Jacka [13], Dubins, Shepp and Shiryaev [9] and Peskir [19]. Problem (20) now models a situation when we incur a running cost and accumulate a reward related to the time spent in a given point as opposed to the reward related to the highest previously visited point. Second, in a similar manner to the results in Peskir ([19], Section 4) for the case of supremum process, our solution to (20) will provide optimal constants in general inequalities involving the local time. We provide a simple illustration of this with Example 3.5. Finally, (20) has a theoretical appeal as an interesting problem which can be solved via a pathwise inequality rather than the more standard approach via a Hamilton–Jacobi–Bellman equation with free boundary.

The approach we use will be based on the representations used in previous sections, where we have made extensive use of the fact that we can construct a local martingale $M_t$ such that $F(L_t) \leq H(B_t) + M_t$. In this section, we can interpret a related martingale as the Snell envelope for the optimal stopping problem; specifically, we are typically able both to provide a meaningful description of the optimal strategy for our problem, and also to write down explicitly the Snell envelope. We believe that being able to get such explicit descriptions of these objects is a strong advantage of this approach.

In this section, we will outline the principle behind this approach, and provide two results, the first of which allows us to provide an upper bound on the problem under very mild conditions on $F$ and $\beta$. The second result gives the value of the problem and an optimal solution under some regularity conditions on $F$ and $\beta$. We then demonstrate through examples that in fact we can find the optimal solution in more general cases. A final example, as stated above, shows how this technique might be used to derive inequalities concerning the local time.

---

[5]Indeed, we can actually solve (20) for $(X_t)$ a regular recurrent time-homogenous diffusion. If $\mathtt{s}$ is the scale function of $(X_t)$ with $s(0) = 0$, then $Y_t = \mathtt{s}(X_t)$ is in natural scale, $\langle Y \rangle_\infty = \infty$ a.s. and the local times at zero of $(X_t)$ and $(Y_t)$ coincide. The problem (20) for $(X_t)$ with a cost function $\beta(x)$ is simply the problem (20) for $(Y_t)$ with the cost function $\beta^Y(y) = \beta(\mathtt{s}^{-1}(y))/(\mathtt{s}'(\mathtt{s}^{-1}(y)))^2$ which can be solved by the methods of the paper.

Writing $H''(x) = 2\beta(x)$, we get

$$H(B_t) = H(0) + \int_0^t H'(B_s)\,dB_s + \int_0^t \beta(B_s)\,ds. \tag{21}$$

In this section, we will further assume that

$$H(x) = \int_0^x \int_0^y 2\beta(z)\,dz\,dy,$$

and therefore $H'$ is continuous and $H(0) = H'(0) = 0$.

Using the results from previous sections, (14) says

$$\begin{aligned}H(B_t) &\geq \gamma_+(L_t)B_t^+ - \gamma_-(L_t)B_t^- + \theta_+(L_t)\mathbf{1}_{B_t \geq 0} + \theta_-(L_t)\mathbf{1}_{B_t < 0} \\ &= M_t^{H,\phi} + \Gamma(L_t) + \theta_+(L_t)\mathbf{1}_{B_t \geq 0} + \theta_-(L_t)\mathbf{1}_{B_t < 0},\end{aligned} \tag{22}$$

where

$$\begin{cases}\gamma_+(l) = H'(\phi_+(l)), & \theta_+(l) = H(\phi_+(l)) - \phi_+(l)H'(\phi_+(l)), \\ \gamma_-(l) = H'(\phi_-(l)), & \theta_-(l) = H(\phi_-(l)) - \phi_-(l)H'(\phi_-(l));\end{cases}$$

and $\Gamma(l) = \int_0^l (\gamma_+(m) - \gamma_-(m))/2\,dm$. In particular,

$$M_t^{H,\phi} = \gamma_+(L_t)B_t^+ - \gamma_-(L_t)B_t^- - \Gamma(L_t)$$

is a local martingale. Combining (21) and (22) we deduce

$$\int_0^t \beta(B_s)\,ds \geq M_t^{H,\phi} - \int_0^t H'(B_s)\,dB_s + \Gamma(L_t) + \theta_+(L_t)\mathbf{1}_{B_t \geq 0} + \theta_-(L_t)\mathbf{1}_{B_t < 0}.$$

Moreover, suppose we can find a solution $(\zeta, \phi_+(\cdot), \phi_-(\cdot))$ to both

$$\begin{aligned}F(l) \leq \zeta &+ \int_0^l \frac{H'(\phi_+(u)) - H'(\phi_-(u))}{2}\,du \\ &+ H(\phi_+(l)) - \phi_+(l)H'(\phi_+(l))\end{aligned} \tag{23}$$

and

$$\begin{aligned}F(l) \leq \zeta &+ \int_0^l \frac{H'(\phi_+(u)) - H'(\phi_-(u))}{2}\,du \\ &+ H(\phi_-(l)) - \phi_-(l)H'(\phi_-(l)).\end{aligned} \tag{24}$$

(Note that, unlike in previous sections, we make no assumption that the functions $\phi_+, \phi_-$ are monotonic.) We can now write

$$F(L_t) - \int_0^t \beta(B_s)\,ds \leq \zeta - M_t^{H,\phi} + \int_0^t H'(B_s)\,dB_s. \tag{25}$$

We note that $N_t^{H,\phi} = \int_0^t H'(B_s)\,dB_s - M_t^{H,\phi}$ is a local martingale with $N_0^{H,\phi} = 0$. In addition, if we define the set

$$\mathcal{T}_\beta = \left\{ \tau : \tau \text{ is a stopping time}, \mathbb{E}\left[\int_0^\tau \beta(B_s)\,ds\right] < \infty \right\}, \tag{26}$$

we deduce from (25) that, when $\tau \in \mathcal{T}_\beta$, $N_{t \wedge \tau}^{H,\phi}$ is bounded below by an integrable random variable, so that it is a supermartingale. Taking expectations, we conclude

$$\mathbb{E}\left[F(L_\tau) - \int_0^\tau \beta(B_s)\,ds\right] \leq \zeta$$

for all stopping times $\tau \in \mathcal{T}_\beta$. In particular we have proved the following result.

PROPOSITION 3.1. *Suppose $F(\cdot)$ and $\beta(\cdot)$ are nonnegative functions, then for any solution $(\zeta, \phi_+(\cdot), \phi_-(\cdot))$ to (23) and (24) we have*

$$\sup_{\tau \in \mathcal{T}_\beta} \mathbb{E}\left[F(L_\tau) - \int_0^\tau \beta(B_s)\,ds\right] \leq \zeta.$$

The arguments which formed the proof of Proposition 3.1 will be important in the sequel. One of our aims will be to obtain an expression for the value of (20) rather than merely a bound. To do this we will need to have equality in (22), (23) and (24), as well as a suitable integrability constraint on the local martingale $N^{H,\phi}$.

If we have equality in (23) and (24) and if we can differentiate suitably, then we must have

$$\phi'_+(l) = \frac{(1/2)(H'(\phi_+(l)) - H'(\phi_-(l))) - F'(l)}{\phi_+(l) H''(\phi_+(l))}, \tag{27}$$

$$\phi'_-(l) = \frac{(1/2)(H'(\phi_+(l)) - H'(\phi_-(l))) - F'(l)}{\phi_-(l) H''(\phi_-(l))}, \tag{28}$$

together with a constraint on the initial values

$$H(\phi_+(0)) - \phi_+(0) H'(\phi_+(0)) = H(\phi_-(0)) - \phi_-(0) H'(\phi_-(0)). \tag{29}$$

Further, equality is attained in (22) exactly on the set where $B_t^+ = \phi_+(L_t)$ or $B_t^- = -\phi_-(L_t)$. Also, since $H(\cdot)$ is convex, the function $H(x) - xH'(x)$ is decreasing in $x$ for $x > 0$, and increasing for $x < 0$. Consequently, we can choose

$$\zeta = F(0) - H(\phi_+(0)) + \phi_+(0)H'(\phi_+(0)) \tag{30}$$

and then also $\zeta = F(0) - H(\phi_-(0)) + \phi_-(0)H'(\phi_-(0))$ where we note that $\zeta$ is increasing if considered as a function of $\phi_+(0)$, and decreasing as a function of $\phi_-(0)$. In particular, we should attempt to minimize $\zeta$ to get a bound which may be attained by the optimal stopping time.

REMARK. Since $\zeta$ is a function of (or determines) our choice of $\phi_-(0)$, $\phi_+(0)$, it seems reasonable to ask how to interpret the solutions of (27)–(29) (i.e., ones with different choices of initial value). In this context there are two possibilities, assuming the relationship in (30) holds.

If we choose $\phi_+(0)$ [and $|\phi_-(0)|$] too small, then we find that the candidate solutions to (27)–(29) hit zero at a finite value $m$, and thereafter the equations no longer make sense. [However, the stopping time $\tau_\phi$ would be optimal for the problem (20) with objective function $F(l \wedge m)$. We will make use of this fact in the sequel.]

Conversely we can ask what happens if we choose initial values for $\phi_+(0)$ and $|\phi_-(0)|$ which are too large. In that case, we can still define a stopping time $\tau_\phi$ but it will not lie in the set $\mathcal{T}_\beta$.

For the main result of this section it will be convenient to introduce the class of finite, positive, continuously differentiable solutions of (27)–(29). For such solutions to exist we will need regularity conditions on $F$ and $\beta$. We define the set

$$\Phi = \left\{ (\phi_+(\cdot), \phi_-(\cdot)) : \begin{array}{c} \phi_+, \phi_- \text{ are solutions of } (27)\text{–}(29), \\ |\phi_+(l)|, |\phi_-(l)| < \infty \quad \forall l \geq 0, \\ |\phi_+(l)|, |\phi_-(l)| > 0 \quad \forall l > 0 \end{array} \right\},$$

where $\phi_+, \phi_-$ are assumed to be continuously differentiable functions on $[0, \infty)$.

THEOREM 3.2. *Suppose $F(\cdot)$ is strictly increasing on $(0, \infty)$ and twice continuously differentiable, and that $\beta$ is positive and continuously differentiable everywhere. Write*

$$V = \sup_{\tau \in \mathcal{T}_\beta} \mathbb{E}\left[ F(L_\tau) - \int_0^\tau \beta(B_s) \, ds \right] \tag{31}$$

*for the value of the optimal stopping problem.*

*Suppose $\Phi = \varnothing$ and*

$$\int_{\mathbb{R}_+} |z|\beta(z) \, dz = \infty = \int_{\mathbb{R}_-} |z|\beta(z) \, dz; \tag{32}$$

*then the value of the problem is infinite.*

*Alternatively, suppose the set $\Phi$ is nonempty; then $\Phi$ contains a minimal [in $\phi_+(0)$] element $(\phi_+(\cdot), \phi_-(\cdot))$, and one of the following is true:*

(i) *The supremum in* (31) *is attained by the stopping time*

$$(33) \qquad \tau = \inf\{t \geq 0 : B_t \notin (\phi_-(L_t), \phi_+(L_t))\}$$

*with corresponding value*

$$(34) \quad V = F(0) + \phi_+(0) H'(\phi_+(0)) - H(\phi_+(0)) = F(0) + \int_0^{\phi_+(0)} z\beta(z)\,dz.$$

(ii) *The stopping time defined in* (33) *is not in $\mathcal{T}_\beta$, but there exists a sequence of stopping times $\tau_N \uparrow \tau$ such that $\tau_N \in \mathcal{T}_\beta$ whose corresponding values*

$$V_N = \mathbb{E}\left[F(L_{\tau_N}) - \int_0^{\tau_N} \beta(B_s)\,ds\right]$$

*increase to $V$, which again is given by* (34).

REMARK. Observe that $(\phi_+(\cdot), \phi_-(\cdot))$, a solution to (27)–(29), is minimal in $\phi_+$ among solutions which do not hit the origin and which remain finite, if and only if it is maximal in $\phi_-$.

Typically in the literature optimal stopping problems like (31) are considered for an arbitrary starting point $(x, l)$ for $(B_t, L_t)$ (e.g., [9, 19]). The fact that a solution is found simultaneously for all possible starting points is a natural consequence, and indeed a necessity, of the fact that the solution method relies heavily on the Markov property, and involves an identification of the value function and the stopping region with a Stefan problem with free boundary. In contrast, in this work we do not rely on Markovian techniques to solve (31) and we consider it only for $(B_t, L_t)$ starting at $(0, 0)$. Naturally once we have the solution given above in Theorem 3.2 the generalization to an arbitrary starting point is straightforward.

Under (32), we can state the conclusions of the theorem in the following way:

*The value $V$ in* (31) *is finite if and only if there exists a minimal solution $(\phi_+, \phi_-)$ to* (27)–(29) *which does not hit the origin, in which case $V$ is given by* (34).

This formulation is parallel to the *maximality principle* described by Peskir [19]. Note also that the solution in (ii) is what Peskir [19] calls an *approximately optimal* solution.

Our aim with this result is not to prove the strongest possible version of Theorem 3.2 since this seems to come at the price of having to account for a variety of idiosyncrasies that these solutions might display. Instead we believe that most, if not all of the restrictions on $F$ and $\beta$ can be weakened in different ways. Examples 3.3 and 3.4 explore this further.

PROOF OF THEOREM 3.2. We first recall classical facts about solutions to ODEs which we use. Under the regularity assumptions on $\beta(\cdot)$ and $F(\cdot)$, the system (27)–(28) can be rewritten as

(35) $$(\phi_-, \phi_+)'(l) = \eta(l, \phi_-(l), \phi_+(l)),$$

with further restriction on the starting point (29), where $\eta(l, x)$ is continuously differentiable on $\mathcal{D} = (0, \infty) \times ((-\infty, 0) \times (0, \infty))$. Classical existence and uniqueness theorems (cf. Arnol'd [2], Section 7.2) imply that for any point $(l, x) \in \mathcal{D}$ there is a unique solution to (35), well defined for some $l_1 < l < l_2$ which goes through $(l, x)$. In particular we can define a solution to (35) for any starting point $x_0$ satisfying (29), and if $\phi_\pm^1$ and $\phi_\pm^2$ are two solutions of (27)–(29) defined on some interval $[0, m)$ with $|\phi_\pm^1(0)| > |\phi_\pm^2(0)|$, then $|\phi_\pm^1(l)| > |\phi_\pm^2(l)|$ for all $l < m$. And in fact, for a fixed $l \in (0, m)$, $(\phi_-^1(l), \phi_+^1(l))$ is a continuous function of the starting point $(\phi_-^1(0), \phi_+^1(0))$ (cf. Arnol'd [2], Section 7.3).

Suppose the set $\Phi$ is nonempty, $(\phi_-^1, \phi_+^1) \in \Phi$. By the above, for any nonzero starting point with $\phi_+(0) \leq \phi_+^1(0)$, and then necessarily $\phi_-(0) \geq \phi_-^1(0)$, the solution to (27)–(29) exists and stays finite and is thus well defined up to the first time $l = m$ when $\phi_-(l)$ or $\phi_+(l)$ hits zero. Note, however, that equalities in (23) and (24) imply

$$H(\phi_+(l)) - \phi_+(l)H'(\phi_+(l)) = H(\phi_-(l)) - \phi_-(l)H'(\phi_-(l))$$

and therefore $\phi_+$ hits zero if and only if $\phi_-$ does. Let $\phi_+(0) = \inf\{\tilde{\phi}_+(0) : (\tilde{\phi}_+, \tilde{\phi}_-) \in \Phi\}$ and $\phi_-(0)$ defined similarly. Naturally $(\phi_-(0), \phi_+(0))$ satisfies (29) and is a starting point of a solution $\phi_\pm(l)$ to (27)–(28). Assume $|\phi_\pm(0)| > 0$. To see that $(\phi_-, \phi_+) \in \Phi$ we have to argue that $\phi_+$ remains strictly positive and $\phi_-$ remains strictly negative. However, the continuity of solutions to (35) in the starting point recalled above implies that $\phi_+$ is the (pointwise) infimum of $\tilde{\phi}_+$ in $\Phi$ and $\phi_-$ is the (pointwise) supremum of $\tilde{\phi}_-$ in $\Phi$. In consequence, the minimal solution $\phi_+$ may only be equal to zero at a point for which $\phi'_+(l) = 0$, and we note that $F'(\cdot) > 0$ rules out this possibility in (27). As noted above, $\phi_-$ can hit zero only if $\phi_+$ does so $\phi_- < 0$.

By the same argument, $\phi_+$ is also the supremum of solutions $\tilde{\phi}_+$ with $\tilde{\phi}_+(0) < \phi_+(0)$ (and then $\phi_-$ is the infimum of $\tilde{\phi}_-$). By definition of $\phi$, such solution $\tilde{\phi}_\pm$ is not in $\Phi$ and thus hits zero in finite time. We now consider these approximating solutions, writing $(\phi_+^m(l), \phi_-^m(l))$ for the solution to (27)–(29) which hits zero at $l = m$. We also write the associated stopping rules $\tau_m$ defined by (33) and $\phi^m$. We argue that $\tau_m$ are optimal for the problems (31) posed for $F(l \wedge m)$. It is clear that $\phi_\pm^m$ are in fact solutions to (27)–(29) with $F(l)$ replaced by $F(l \wedge m)$. In view of the arguments which led to Proposition 3.1, the only property we need to demonstrate is that the supermartingale $N_{t \wedge \tau_m}^{H,\phi} = \int_0^{t \wedge \tau_m} H'(B_s) \, dB_s - M_{t \wedge \tau_m}^{H,\phi}$ is in fact a UI martingale. The important point to note here is that, by construction, the stopping times $\tau_m$ are smaller than $\inf\{t \geq 0 : L_t \vee |B_t| \geq J\}$ for some $J$. It follows immediately that $M_{t \wedge \tau_m}^{H,\phi}$ is a UI martingale. On the other hand the local martingale $\int_0^{t \wedge \tau_m} H'(B_s) \, dB_s$ has quadratic variation which is bounded by

$$\int_0^{\tau_m} H'(B_s)^2 \, ds \leq K \tau_m,$$

for some $K > 0$. We know that $B_{t \wedge \tau_m}$ satisfies the conditions of the Azéma–Gundy–Yor theorem [3], Theorem 1b, since $\tau_m$ is bounded by the hitting time of $\{-J, J\}$. As a consequence, $\int_0^{t \wedge \tau_m} H'(B_s) \, dB_s$ also satisfies the conditions of the Azéma–Gundy–Yor theorem and thus is a uniformly integrable martingale.

This procedure results in a sequence of stopping times, optimal for the problems posed with gain function $F(l \wedge m)$, and with values increasing to

$$V = F(0) - H(\phi_+(0)) + \phi_+(0) H'(\phi_+(0))$$

from which we deduce that we are in either case (i) or (ii).

We must also consider the case where the minimal solution is a solution with $\phi_+(0) = \phi_-(0) = 0$, and corresponding stopping time $\tau \equiv 0$. It is then trivial to apply Proposition 3.1 to this solution to deduce that this is the optimal solution, with $V = F(0)$.

It remains to prove the initial statement of the theorem. Let $(\phi_-(\cdot), \phi_+(\cdot))$ be a solution of (27)–(29) with some nonzero starting point. Since $\Phi = \varnothing$ this solution has to either hit zero or explode in finite time. We show that the latter is impossible. Define

$$G(z) = 2 \int_{\phi_-(z)}^{\phi_+(z)} |u| \beta(u) \, du.$$

We have

$$\begin{aligned} G'(z) &= 2 \phi_+'(z) \phi_+(z) \beta(\phi_+(z)) + 2 \phi_-'(z) \phi_-(z) \beta(\phi_-(z)) \\ &= H'(\phi_+(z)) - H'(\phi_-(z)) - 2 F'(z) \\ &\leq 2 \int_{\phi_-(z)}^{\phi_+(z)} \beta(u) \, du \\ &\leq 2 \int_{\phi_-(0)}^{\phi_+(0)} \beta(u) \, du + \frac{2}{\phi_+(0)} \int_0^{\phi_+(z)} u \beta(u) \, du \end{aligned}$$

$$+ \frac{2}{|\phi_-(0)|} \int_{\phi_-(z)}^{0} |u|\beta(u)\,du$$
$$\leq C_1 + C_2 G(z),$$

where $C_1 = 2\int_{\phi_-(0)}^{\phi_+(0)} \beta(u)\,du$ and $C_2^{-1} = \min\{\phi_+(0), |\phi_-(0)|\}$. It follows from Gronwall's lemma that $G(u) \leq C_1 e^{C_2 u}$; combining this with (32) we see that neither $\phi_+$ nor $\phi_-$ can explode. We conclude that $\phi_\pm$ are bounded and hit zero in finite time. In consequence, the stopping time associated via (33) is in $\mathcal{T}_\beta$. We can apply Proposition 3.1 to see that

$$V \geq F(0) + \lim_{\phi_+(0)\uparrow\infty} [\phi_+(0)H'(\phi_+(0)) - H(\phi_+(0))]$$

and the value of the problem is infinite since (32) ensures $\int_0^x z\beta(z)\,dz = xH'(x) - H(x)$ increases to infinity. $\square$

We finish this section with three examples. The first and second examples aim to demonstrate that, although we require relatively strong conditions in order to apply Theorem 3.2, the techniques and principles of the result are much more generally applicable. The first example also shows that the condition (32) cannot be weakened in general. A final example connects the results with established inequalities concerning the local time.

EXAMPLE 3.3. The initial example considers the case where $F(l) = l$ and $\beta(x) > 0$ is continuous. In this setting, there are three possible types of behavior, depending on the value of

$$c = \int_\mathbb{R} \beta(x)\,dx.$$

For $c < 1$, every solution of (27) [resp. (28)] will have a strictly negative (resp. positive) gradient, and will therefore hit the origin in finite time. Consider $\tau = \inf\{t \geq 0 : L_t = 1\}$ and $\tau_N = \inf\{t \geq 0 : |B_t| \geq N \text{ or } L_t = 1\}$. Then $\mathbb{P}(B_{\tau_N} = 0) = \exp(-1/N)$, and by applying Itô's formula to $H(B_t)$ and using monotone convergence, and using the fact that $H'$ is increasing and bounded,

$$\mathbb{E}\left[\int_0^\tau \beta(B_s)\,ds\right] = \lim_{N\to\infty} \mathbb{E}H(B_{\tau_N})$$
$$= \lim_{N\to\infty}\left[\frac{H(N) + H(-N)}{2}\right](1 - e^{-1/N})$$
$$= \lim_{N\to\infty}\left(\frac{H'(N) - H'(-N)}{2}\right)$$
$$= c.$$

As a consequence, there is a positive, finite gain of $1-c$ from simply running until the local time reaches 1. This process can then be continued, waiting until the local time reaches an arbitrary level, giving an infinite value to the problem.

The interesting case occurs when $c=1$. In this setting, we see that the above argument fails—the strategy of waiting until the local time reaches 1 has no average gain. However, there is still value to the problem; we can apply Proposition 3.1 (when we interpret $\phi_+ = \phi_- = \infty$ suitably) to deduce that the value of the problem is at most

$$(36) \qquad \max\left(\lim_{x\to\infty}[xH'(x) - H(x)], \lim_{x\to-\infty}[xH'(x) - H(x)]\right),$$

and we note that this expression can be infinite. Now consider the payoff from running to exit of the interval $(-\alpha N, N)$ for some $\alpha \in (0, \infty)$. From the martingale property of $L_t - |B_t|$ this is easily seen to be

$$\frac{2\alpha N}{1+\alpha} - H(N)\frac{\alpha}{1+\alpha} - H(-\alpha N)\frac{1}{1+\alpha}.$$

Using the fact that $1 = \frac{1}{2}(H'(\infty) - H'(-\infty)) \geq \frac{1}{2}(H'(N) - H'(-\alpha N))$, we can bound the last expression from below by

$$\frac{\alpha}{1+\alpha}[NH'(N) - H(N)] + \frac{1}{1+\alpha}[(-\alpha N)H'(-\alpha N) - H(-\alpha N)].$$

Since $\alpha$ was arbitrary, we conclude that the solutions for sufficiently large $N$ and sufficiently large/small $\alpha$ are approximately optimal and the value in (36) is obtained in the limit.

Finally, we consider the case $c > 1$. We want to find a pair $\phi_-, \phi_+$ such that $\frac{1}{2}(H'(\phi_+) - H'(\phi_-)) = 1$ and such that (29) holds. We can rewrite this in terms of $\beta$:

$$(37) \qquad \int_{\phi_-}^{\phi_+} s\beta(s)\,ds = 0 \quad \text{subject to} \quad \int_{\phi_-}^{\phi_+} \beta(s)\,ds = 1.$$

If this has a finite solution, then we can apply Theorem 3.2 to conclude that the value function is given by $\phi_+ H'(\phi_+) - H(\phi_+) = \phi_- H'(\phi_-) - H(\phi_-)$ and the optimal strategy is to stop on exit from a finite interval.

Otherwise (37) has no solution with both $\phi_+$ and $\phi_-$ finite and either

$$(38) \qquad \int_{-\infty}^{\phi_+} s\beta(s)\,ds \geq 0 \qquad \text{for } \phi_+ \text{ the solution of } \int_{-\infty}^{\phi_+} \beta(s)\,ds = 1,$$

or

$$(39) \qquad \int_{\phi_-}^{\infty} s\beta(s)\,ds \leq 0 \qquad \text{for } \phi_- \text{ the solution of } \int_{\phi_-}^{\infty} \beta(s)\,ds = 1.$$

For the former case we must have

(40) $$\lim_{x\downarrow-\infty} xH'(x) - H(x) < \lim_{x\uparrow\infty} xH'(x) - H(x)$$

whereas for (39) we must have the reverse. [Note that if both sides of (40) are infinite, then we must have a finite solution to (37).] Assuming (38) holds, it follows from Proposition 3.1 that $\phi_+ H'(\phi_+) - H(\phi_+)$ is an upper bound on the value function; that this bound can be attained in the limit follows from consideration of stopping rules which are the first exit times from intervals of the form $(-N, \phi_+)$. If the inequality in (40) is reversed, then the value function is $\phi_- H'(\phi_-) - H(\phi_-)$ where $\phi_-$ solves the integral equality in (39).

Using this setup one can easily construct an example when the value is finite even though $\int_\mathbb{R} |z|\beta(z)\,dz = \infty$ which shows that both $\int_{\mathbb{R}_-} |z|\beta(z)\,dz = \infty$ and $\int_{\mathbb{R}_+} z\beta(z)\,dz = \infty$ are needed in general to ensure the last statement of Theorem 3.2.

EXAMPLE 3.4. The main aim of this example is to demonstrate that the above ideas can lead to meaningful solutions to optimal stopping problems even if the functions $F$ and $H$ and the resulting minimal solution to (27)–(29) are not "nice." In particular, we shall give an example where the minimal solution is finite only on a bounded interval which excludes the origin, and where $F$ is not nondecreasing, but where the value of the problem is finite.

Specifically, we consider the optimal stopping problem (31) where the (slightly contrived) functions $F$ and $\beta$ are defined by

$$F'(l) = \begin{cases} (3-2l)e^{-1/2l^2(l-2)^2}, & l < 2, \\ 1, & l \geq 2, \end{cases}$$

with $F(0) = 0$ and

$$2\beta(x) = H''(x) = |x|^{-3} e^{-1/(2x^2)}.$$

We also obtain

$$H'(x) = e^{-1/(2x^2)}, \qquad H(x) = xe^{-1/(2x^2)} - \int_{1/x}^{\infty} e^{-1/2z^2}\,dz.$$

Noting that the problem is symmetric (and therefore dropping the subscripts to denote positive and negative solutions), (27) becomes

(41) $$\phi'(l) = \begin{cases} \phi(l)^2 \left(1 - (3-2l)\exp\left\{\dfrac{1}{2\phi(l)^2} - \dfrac{1}{2}l^2(2-l)^2\right\}\right), & l < 2, \\ \phi(l)^2 \left(1 - \exp\left\{\dfrac{1}{2\phi(l)^2}\right\}\right), & l \geq 2. \end{cases}$$

Define the function $\phi_0(l)$ via

$$\phi_0(l) = \begin{cases} \dfrac{1}{l(2-l)}, & l \in (0,2), \\ \infty, & \text{otherwise;} \end{cases}$$

then $\phi_0$ is a solution of (41) for $l < 2$. Moreover if for $l > 2$ we use the natural definitions $H'(\phi) = 1$ and $H(\phi) - \phi H'(\phi) = -\sqrt{\pi/2}$ when $\phi = \infty$, then

$$F(l) = \int_0^l H'(\phi_0(u))\,du + H(\phi_0(l)) - \phi_0(l)H'(\phi_0(l)) + \sqrt{\frac{\pi}{2}}$$

for all $l \geq 0$.

We will show that the stopping time $\tau = \inf\{t \geq 0 : |B_t| \geq \phi_0(L_t)\}$ is *approximately optimal* in the sense described in the remarks below Theorem 3.2. Specifically, we consider a sequence of solutions $\phi_m$ increasing to $\phi_0$ which have expected values increasing to $\sqrt{\frac{\pi}{2}}$, and show further that no solution can improve on this bound. The solutions $\phi_m(l)$ and $\phi_0(l)$ are shown in Figure 3.

It is straightforward to check that the solutions $\phi_m$ do indeed increase to $\phi_0$, and that each $\phi_m$ is a well-defined solution to (31) where the function $F$ is replaced by $F(l) \wedge m'$ for some $m'$. The resulting sequence provides an increasing set of stopping times in $\mathcal{T}_\beta$, with values increasing to $\sqrt{\frac{\pi}{2}}$.

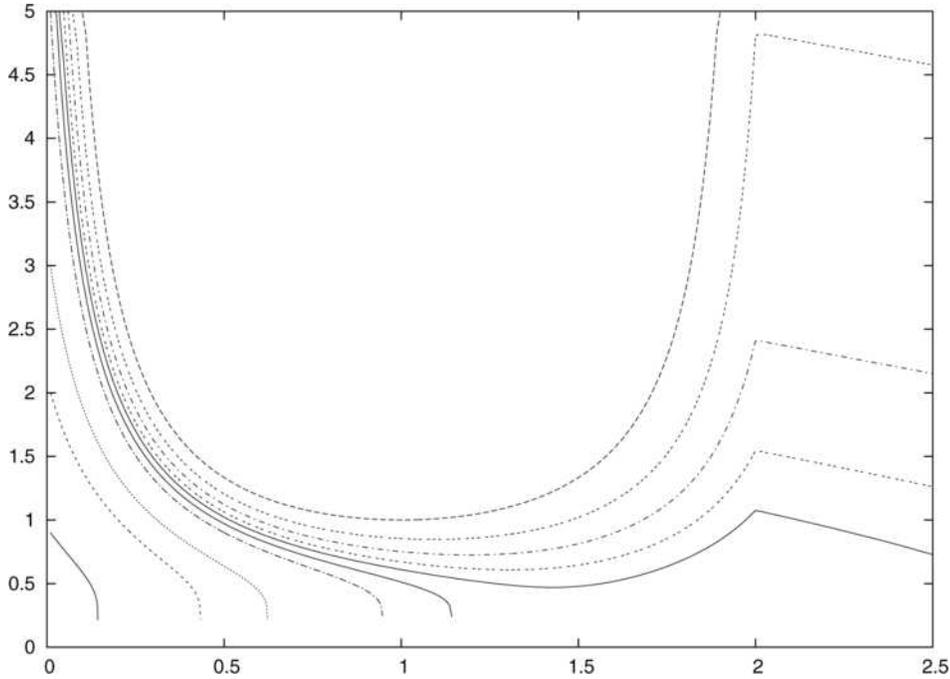

FIG. 3. *Solutions to* (41) *with different initial values. The top curve shows the function* $\phi_0(l)$.

To conclude that this does give the supremum, we suppose for a contradiction that there exists $\tau \in \mathcal{T}_\beta$ with larger expected value. Then for sufficiently large $m'$ we have

$$\mathbb{E}\left[F(L_\tau) \wedge m' - \int_0^\tau \beta(B_s)\, ds\right] > \sqrt{\frac{\pi}{2}},$$

contradicting the optimality of the $\tau_m$.

Finally, we note that the value of the problem for $F(l) = (l - l_0)_+$ is $(\sqrt{\frac{\pi}{2}} - l_0)_+$, so that $F(l) = l$ has the same value as our original problem, while $F(l) = (l - \sqrt{\frac{\pi}{2}})_+$ has zero value. The value for $F(l)$ actually follows from Example 3.3 since $H'(\infty) = -H'(-\infty) = 1$ and this also provides additional intuition behind the above results.

EXAMPLE 3.5. We end by demonstrating how our techniques may be used to recover the well-known inequality

(42) $$\mathbb{E}L_\tau^p \leq p^p \mathbb{E}|B_\tau|^p, \qquad p > 1,$$

valid for all stopping times $\tau$ such that the $\mathbb{E}|B_\tau|^p = \frac{p(p-1)}{2}\mathbb{E}\int_0^\tau |B_s|^{p-2}\, ds$. Fix $p > 1$ and consider $F(l) = \frac{l^p}{p}$ and $\beta_c(x) = c|x|^{p-2}/2$ for some $c > 0$. The function $H_c(x) = \frac{c}{p(p-1)}|x|^p$ is symmetric so that $\phi_- = \phi_+$. One can easily verify that $\phi(x) = ax$ satisfies (27) if and only if $ca^p - \frac{c}{p-1}a^{p-1} + 1 = 0$ and that this equation has two solutions only for $c \geq c_{\min} = p^p(p-1)$. As $\phi$ is linear $\int_{0+} ds/\phi(s) = \infty$ and the resulting stopping time $\tau_V = 0$ a.s. Thus the value $V$ of the optimal stopping problem (31) associated with $F$ and $\beta_c$ is zero. Consequently

(43) $$\mathbb{E}\frac{L_\tau^p}{p} \leq \frac{c}{p(p-1)}\mathbb{E}|B_\tau|^p, \qquad \tau \in \mathcal{T}_\beta, c \geq c_{\min}$$

and we recover (42) on taking $c = c_{\min}$.

**Acknowledgment.** We thank two anonymous referees whose comments helped us improve the first version of the paper.

A. M. G. Cox  
Department of Mathematical Sciences  
University of Bath  
Bath BA2 7AY  
United Kingdom  
E-mail: A.M.G.Cox@bath.ac.uk

D. Hobson  
Department of Statistics  
University of Warwick  
Coventry CV4 7AL  
United Kingdom  
E-mail: D.Hobson@warwick.ac.uk

J. Obłój  
Department of Mathematics  
Imperial College London  
London SW7 2AZ  
United Kingdom  
E-mail: j.obloj@imperial.ac.uk